\documentclass[twoside,11pt]{article}

\usepackage{amsmath,amssymb, amsthm}
\usepackage{graphicx}

\newtheorem{thm}{Theorem}

\newtheorem{lem}{Lemma}
\newtheorem{prop}{Proposition}
\theoremstyle{definition}
\newtheorem{defn}{Definition}
\theoremstyle{remark}
\newtheorem{rem}{Remark}
\numberwithin{equation}{section}
\begin{document}
	\pagestyle{myheadings} \markboth{ \rm \centerline {Mateusz Kubiak and Bogdan Szal}} 
	{\rm \centerline { }}
	
	\begin{titlepage}
		\title{\bf {A sufficient condition for uniform convergence of trigonometric series with $p$-bounded variation coefficients}}
		\author {Mateusz Kubiak and Bogdan Szal \\
			{\small University of Zielona G\'{o}ra,}\\
			{\small Faculty of Mathematics, Computer Science and Econometrics,}\\
			{\small 65-516 Zielona G\'{o}ra, ul. Szafrana 4a, Poland} \\ 
			{\small e-mail: M.Kubiak@wmie.uz.zgora.pl, B.Szal @wmie.uz.zgora.pl}}
	\end{titlepage}
	
	\date{}
	\maketitle
	\begin{abstract}
		In this paper we consider trigonometric series with $p-$bounded variation coefficients. We presented a sufficient condition for uniform convergance of such series in case $p>1$. This condition is significantly weaker than these obtained in the results on this subject known in the literature.
	\end{abstract}
		
\section{Introduction}
	It is well known that there is a great number of interesting results in
	Fourier analysis established by assuming monotonicity of Fourier coefficients. The
	following classical convergence result can be found in many monographs (see for example \cite{Chaundy} and \cite{Zygmund}).
	
	\begin{thm}\label{pierwsze twierdzenie}
		Suppose that $b_{n}\geq b_{n+1}$ and $b_{n}\rightarrow 0$ as $n \rightarrow \infty$. Then a necessary
		and sufficient condition for the uniform convergence of the series%
		\begin{equation}
		\sum\limits_{n=1}^{\infty }b_{n}\sin nx  \label{1}
		\end{equation}%
		is $nb_{n}\rightarrow 0$ as $n \rightarrow \infty$.
	\end{thm}

	This result has been generalized by weakening the monotonicity conditions of
	the coefficient. We present below a historical outline of the generalizations of this theorem.  
	
	In 2001 Leindler defined (see \cite{Leindler1} and \cite{Leindler_new}) a new class of sequences named as
	sequences of Rest Bounded Variation, briefly denoted by $RBVS$, i.e.,%
	\begin{equation*}
	RBVS=\left\{ a:=\left( a_{n}\right) \in 
	\mathbb{C}
	:\sum\limits_{n=m}^{\infty }\left\vert \Delta_1a_n\right\vert \leq
	C \left\vert a_{m}\right\vert \text{ for all }m\in 
	\mathbb{N}
	\right\} ,
	\end{equation*}%
	where here and throughout the paper $C=C\left( a\right) $ always indicates a
	constant only depending on $a$ and $\Delta_ra_n = a_n-a_{n+r}$ for $r \in \mathbb{N}$.
	
	Denote by $MS$ the class of monotone decreasing sequences, then it is clear  that%
	\begin{equation}\label{1.2}
	MS\subsetneq RBVS. 
	\end{equation}%
	Further,  Tikhonov introduced a class of General Monotone Sequences $GMS$ defined as
	follows (see \cite{Tikhonov3}):%
	\begin{equation*}
	GMS=\left\{ a:=\left( a_{n}\right) \in 
	\mathbb{C}
	:\sum\limits_{n=m}^{2m-1}\left\vert\Delta_1a_n\right\vert \leq C \left\vert a_{m}\right\vert \text{ for all }m\in 
	\mathbb{N}
	\right\} .
	\end{equation*}%
	It is clear that%
	\begin{equation}\label{1.3}
	RBVS\subsetneq GMS\text{.}
	\end{equation}

	The class of $GMS$ was generalized by Tikhonov (see \cite{Tikhonov2}) and independently by Zhou, Zhou and Yu (see \cite{zhou}) to the class of Mean Value Bounded Variation Sequences ($MVBVS$). We say that a sequence $a:=\left( a_{n}\right) $ of complex numbers is said to be $MVBVS$ if there exists $\lambda\geq 2$ such that
	\begin{equation*}
	\sum\limits_{k=n}^{2n-1}\left\vert \Delta_1a_k\right\vert \leq
	\frac{C}{n} \sum\limits_{k=[n/\lambda]}^{\lambda n} |a_{k}|
	\end{equation*}%
	holds for $n \in \mathbb{N}$, where $[x]$ is the integer part of $x$. They proved also in \cite{zhou} that  
	\begin{equation}\label{1.4}
	GMS \subsetneq MVBVS.
	\end{equation}
	
	Theorem 1 was generalized for the class $RBVS$ in \cite{Leindler1}, for the class $GMS$ in \cite{Tikhonov3} and for the class $MVBVS$ in \cite{zhou}.
	
	Next, Tikhonov (\cite{Tikhonov1}, \cite{Tikhonov2}, \cite{Tikhonov3}) and Leindler (\cite{Leindler2}) defined the class of $\beta\text{ -}$ general monotone
	sequences as follows:
	
	\begin{defn}
		Let $\beta :=\left( \beta _{n}\right) $ be a nonnegative sequence. The
		sequence of complex numbers $a:=\left( a_{n}\right) $ is said to be $\beta\text{ -}$ 
		general monotone, or $a\in GM\left( \beta \right) $, if the relation%
		\begin{equation*}
		\sum\limits_{n=m}^{2m-1}\left\vert \Delta_1a_n\right\vert \leq C \beta _{m}
		\end{equation*}%
		holds for all $m \in \mathbb{N}$.
	\end{defn}
	
	In the paper \cite{Tikhonov2} Tikhonov considered i.e. the following examples of the
	sequences $\beta _{n}:$
	
	(1) $_{1}\beta _{n}=\left\vert a_{n}\right\vert ,$
	
	(2) $_{2}\beta _{n}=\sum\limits_{k=\left[ n/c%
		\right] }^{\left[ cn\right] }\frac{\left\vert a_{k}\right\vert }{k}$ for
	some $c>1$.
	
	It is clear that $GM\left( _{1}\beta \right) =GMS$. Moreover, Tikhonov showed in \cite{Tikhonov2} that
	\begin{equation}
	GM\left( _{1}\beta \right) \subsetneq GM\left( _{2}\beta \right) \equiv \text{MVBVS}.  \label{1c}
	\end{equation}
	
	Tikhonov proved also in \cite{Tikhonov2} the following theorem:
	\begin{thm}
		Let a sequence $(b_n) \in GM(\text{}_2\beta)$. If $n|b_n| \rightarrow 0$ as $n\rightarrow \infty$, then the series (1.1) converges uniformly.
	\end{thm}
	
	Further, Szal defined a new class of sequences in the following way (see \cite{Szal}):
	
	\begin{defn}
		Let $\beta :=\left( \beta _{n}\right) $ be a nonnegative sequence and $r$ a
		natural number. The sequence of complex numbers $a:=\left( a_{n}\right) $ is
		said to be $\left( \beta ,r\right) -$general monotone, or $a\in GM\left(
		\beta ,r\right) $, if the relation%
		\begin{equation*}
		\sum\limits_{n=m}^{2m-1}\left\vert \Delta_ra_n\right\vert \leq C \beta _{m}
		\end{equation*}%
		holds for all $m \in \mathbb{N}$.
	\end{defn}
	
	It is clear that $GM\left( \beta ,1\right) \equiv GM\left( \beta \right) $. Moreover, it is easy to show that the sequence
	\begin{equation*}
	a_n = \frac{\left(-1\right)^n}{n}
	\end{equation*}
	belongs to $GM\left( {}_1\beta ,2\right)$ and does not belong to $GM\left( {}_1\beta\right)$. This example shows that the class $GM\left( {}_1\beta\right)$ is essentially wider than the class $GM\left( {}_1\beta\right)$. In \cite{Szal} Szal showed more general relations
	\begin{equation}\label{1.61}
	GM({}_2\beta,1) \subsetneq GM({}_2\beta,r)
	\end{equation}
	for all $r >1$.
	
	In the paper \cite{Szal} Szal generalized Theorem \ref{pierwsze twierdzenie} by proving the following theorem.
	
	\begin{thm} \text{(\cite{Szal})}
		Let a sequence $\left( b_{n}\right) \in GM\left( {}_2\beta,r\right) $, where $r \in \mathbb{N}$. If $n|b_n| \rightarrow 0$ as $n \rightarrow \infty$ and 
		\begin{equation*}
		\sum\limits_{n=1}^{\infty}\sum\limits_{k=1}^{[r/2]}\left\vert b_{r\cdot n +k} - b_{r\cdot n + r - k} \right\vert < \infty \text{ for } r\geq 3,
		\end{equation*}
		then the series (\ref{1}) converges uniformly.
	\end{thm}
	
	In the paper \cite{Korus_SBVS_2} Kórus defined a new class of sequences in the following way:
	\begin{defn}
		The sequence of complex numbers $a :=\left( a _{n}\right)$  is in the class $SBVS_2$ (Supremum Bounded Variation Sequence), if the relation
		\begin{equation*}
		\sum\limits_{n=m}^{2m-1}\left\vert \Delta_1 a_n\right\vert \leq C \frac{1}{n} \underset{m\geq b(n)}{\sup} \sum\limits_{k=m}^{2m} |a_{k}|
		\end{equation*}%
		holds for all $m\in \mathbb{N}$, where $(b(n))$ is a nonegative sequence tending monotonically to infinity depending only on $a$.
	\end{defn}
	In the paper \cite{Korus_SBVS_2} Kórus also proved the following therorem:
	\begin{thm}
		Let a sequence $(b_n) \in SBVS_2$. If $n|b_n| \rightarrow 0$ as $n\rightarrow \infty$, then the series (\ref{1}) converges uniformly. 
	\end{thm}

	Next Tikhonov and Liflyand defined a class of $GMS_p(\beta)$ in the following way (see \cite{Liflyand_Tikhonov}, \cite{Tikhonov1_new}):

	\begin{defn}
		Let  $\beta = (\beta_n)$ be a nonnegative sequence and $p$ a positive real number. We say that a sequence of complex numbers $a=(a_n) \in GMS_p(\beta)$ if the relation
		\begin{equation*}
		(\sum\limits_{n=m}^{2m-1}\left\vert \Delta_1a_n\right\vert^p)^{\frac{1}{p}} \leq K(a)  \beta_{m}
		\end{equation*}
		holds for all $m\in \mathbb{N}$.
	\end{defn}
	It is clear that $GMS_1(\beta)=GM(\beta)$.
	
	The latest class of sequences was defined by Kubiak and Szal in \cite{Kubiak} as follows:
	\begin{defn}
		Let $\beta :=\left( \beta _{n}\right) $ be a nonnegative sequence, $r$ a
		natural number and $p$ a positive real number. The sequence of complex numbers $a:=\left( a_{n}\right) $ is
		said to be $\left( p,\beta ,r\right) -$ general monotone, or $a\in GM\left(p,
		\beta ,r\right) $, if the relation%
		\begin{equation*}
		(\sum\limits_{n=m}^{2m-1}\left\vert \Delta_ra_n\right\vert^p)^{\frac{1}{p}} \leq C \beta_m
		\end{equation*}%
		holds for all $m\in \mathbb{N}$.
	\end{defn}
	
	It is clear that $GM(p,\beta,1)=GMS_p(\beta)$ and $GM(1,\beta,r) = GM(\beta,r)$.

	Further we will consider the following sequence:
	\begin{equation*}
	{}_3\beta_n(q)=\frac{1}{n} \underset{m\geq b(n)}{\sup}m \left(\frac{1}{m} \sum\limits_{k=m}^{2m} |a_{k}|^q   \right)^{\frac{1}{q}},
	\end{equation*}
	where $(a_n) \subset \mathbb{C}, a_n\rightarrow 0$ as $n\rightarrow \infty$, $q> 0, \lambda \geq 2$, $(b(n))$ is a nonnegative sequence such that $b(n) \nearrow $ and $b(n) \rightarrow \infty$ as $n \rightarrow \infty$. It is clear that $SBVS_2 = GM(1, {}_3\beta(1),1)$.
	
	In the further part of our paper we will consider the following series:
	\begin{equation}\label{sinus}
	\sum\limits_{n=1}^{\infty} b_n \sin (cnx),
	\end{equation}
	
	\begin{equation}\label{cosinus}
	\sum\limits_{n=1}^{\infty} a_n \cos (cnx),
	\end{equation}
	\begin{equation}\label{exp}
	\sum\limits_{n=1}^{\infty} c_n e^{ickx},
	\end{equation}
	where $c >0.$

	In the paper \cite{Kubiak} Kubiak and Szal showed the following embedding relations:
	\begin{thm}
		Let $q >0, r\in \mathbb{N}$ and $0 < p_1 \leq p_2$. Then 
		\begin{equation*}
		GM(p_1,{}_3\beta(q),r) \subseteq GM(p_2,{}_3\beta(q),r).
		\end{equation*}
	\end{thm}
	
	\begin{thm}
		Let $p \geq 1$, $q > 0$, $r_1,r_2 \in \mathbb{N}$, $r_1 \leq r_2$. If $r_1 | r_2$, then 
		\begin{equation*}
		GM(p,{}_3\beta(q),r_1) \subseteq GM(p,{}_3\beta(q),r_2).
		\end{equation*}
	\end{thm}

	Moreover they proved in \cite{Kubiak} the following generalization of Theorem \ref{pierwsze twierdzenie}:
	
	\begin{thm}\label{twierdzenie_Kubiak_Szal_sinus}
		Let a sequence $\left( b_{n}\right) \in GM\left(p, {}_{3}\beta(q) ,r\right) $,
		where $p,q\geq 1$, $r \in \mathbb{N}$ and $b(n)\geq n$ for $n \in \mathbb{N}$. If 
		
		\begin{equation}\label{prof 12}
		n^{2-\frac{1}{p}}\left\vert b_{n}\right\vert \rightarrow 0\text{ as } %
		n\rightarrow \infty 
		\end{equation}
		and 
		\begin{equation}\label{2.4}
		\sum\limits_{k=1}^{\infty} b_k \sin \left(\frac{2l \pi }{r} k\right) <\infty ,\text{ for }r\geq 3,
		\end{equation}%
		for all $l= 1,...,[\frac{r}{2}]-1 $ when $r$ is an even number and $l= 1,...,[\frac{r}{2}] $ when $r$ is an odd number, then the series \eqref{sinus} is uniformly convergent.
	\end{thm}
	
	\begin{thm}\label{twierdzenie_Kubiak_Szal_cosinus}
		Let a sequence $\left( a_{n}\right) \in GM\left(p, {}_{3}\beta(q) ,r\right) $,
		where $p,q\geq 1$, $r \in \mathbb{N}$ and $b(n)\geq n$ for $n \in \mathbb{N}$. If 
		\begin{equation*}
		n^{2-\frac{1}{p}}\left\vert a_{n}\right\vert \rightarrow 0\text{ as	}n\rightarrow \infty
		\end{equation*} 
		and 
		\begin{equation}\label{2.41}
		\sum\limits_{k=1}^{\infty} a_k \cos \left(\frac{2l \pi }{r} k\right) <\infty ,
		\end{equation}%
		for all $l= 0,1,...,[\frac{r}{2}] $, then the series \eqref{cosinus} is uniformly convergent.
	\end{thm}
	
	\begin{thm}\label{twierdzenie_Kubiak_Szal_exp}
		Let a sequence $\left( c_{n}\right) \in GM\left(p, {}_{3}\beta(q) ,r\right) $,
		where $p,q\geq 1$, $r \in \mathbb{N}$ and $b(n)\geq n$ for $n \in \mathbb{N}$. If 
		\begin{equation*}
		n^{2-\frac{1}{p}}\left\vert c_{n}\right\vert \rightarrow 0\text{ as }n\rightarrow \infty 
		\end{equation*} 
		and 
		\begin{equation}\label{2.42}
		\sum\limits_{k=1}^{\infty} c_k e^{\left(\frac{2l \pi }{r} k\right)i} <\infty ,
		\end{equation}%
		for all $l= 0,1,...,[\frac{r}{2}] $, then the series \eqref{exp} is uniformly convergent.
	\end{thm}

	In this paper we will show that Theorem \ref{twierdzenie_Kubiak_Szal_sinus}, \ref{twierdzenie_Kubiak_Szal_cosinus}, \ref{twierdzenie_Kubiak_Szal_exp} are true under weakened assumptions in case $p>1$.
	
\section{Main results}
	We have the following results:

	\begin{thm}\label{twierdzenie_sinus_nowe}
		Let a sequence $\left( b_{n}\right) \in GM\left(p, {}_{3}\beta(q) ,r\right) $,
		where $q\geq 1$, $p>1$ and $r \in \mathbb{N}$. If 
		
		\begin{equation}\label{warunek_n_ln_n}
		n \ln n\left\vert b_{n}\right\vert \rightarrow 0\text{ as } %
		n\rightarrow \infty 
		\end{equation}
		and 
		\begin{equation}\label{zbiezne_w_punkach}
		\sum\limits_{k=1}^{\infty} b_k \sin \left(\frac{2l \pi }{r} k\right) <\infty ,\text{ for }r\geq 3,
		\end{equation}%
		for all $l= 1,...,[\frac{r}{2}]-1 $ when $r$ is an even number and $l= 1,...,[\frac{r}{2}] $ when $r$ is an odd number, then the series \eqref{sinus} is uniformly convergent.
	\end{thm}

	\begin{prop}\label{uwaga_przyklad}
		There exist an $x_0 \in \mathbb{R}$ and a sequence $(b_n) \in GM(p, {}_3\beta(1), 3)$ for $p>1$ with the properties $nb_n \rightarrow 0 $ as $n \rightarrow \infty$ and $(b_n) \notin GM(1, {}_3\beta(1),3)$, for which the series (\ref{sinus}) is divergent in $x_0$.
	\end{prop}

	\begin{thm}\label{twierdzenie_cosinus_nowe}
		Let a sequence $\left( a_{n}\right) \in GM\left(p, {}_{3}\beta(q) ,r\right) $,
		where $q\geq 1$, $p>1$ and $r \in \mathbb{N}$. If 
		\begin{equation*}
		n \ln n\left\vert a_{n}\right\vert \rightarrow 0\text{ as } %
		n\rightarrow \infty
		\end{equation*}
		and 
		\begin{equation}\label{2.41}
		\sum\limits_{k=1}^{\infty} a_k \cos \left(\frac{2l \pi }{r} k\right) <\infty ,
		\end{equation}%
		for all $l= 0,1,...,[\frac{r}{2}] $, then the series \eqref{cosinus} is uniformly convergent.
	\end{thm}
	
	\begin{thm}\label{twierdzenie_exp_nowe}
		Let a sequence $\left( c_{n}\right) \in GM\left(p, {}_{3}\beta(q) ,r\right) $,
		where $q\geq 1$, $p>1$ and $r \in \mathbb{N}$. If 
		\begin{equation*}
		n \ln n\left\vert c_{n}\right\vert \rightarrow 0\text{ as } %
		n\rightarrow \infty
		\end{equation*}
		and 
		\begin{equation}\label{2.42}
		\sum\limits_{k=1}^{\infty} c_k e^{\left(\frac{2l \pi }{r} k\right)i} <\infty ,
		\end{equation}%
		for all $l= 0,1,...,[\frac{r}{2}] $, then the series \eqref{exp} is uniformly convergent.
	\end{thm}

	\begin{rem}
		It is clear that if a sequence $(b_n)$ satisfies the condition (\ref{warunek_n_ln_n}) then it fulfills the condition (\ref{prof 12}) with $p>1$, too. Therefore, from Theorem \ref{twierdzenie_sinus_nowe} we get Theorem \ref{twierdzenie_Kubiak_Szal_sinus} is case $p>1$. The same remark applies to Theorems \ref{twierdzenie_cosinus_nowe}, \ref{twierdzenie_Kubiak_Szal_cosinus} and Theorems \ref{twierdzenie_exp_nowe}, \ref{twierdzenie_Kubiak_Szal_exp}, respectively.
	\end{rem}
	
\section{Lemma}
	
	Denote, for $r \in \mathbb{N}$ and $k=0,1,2...$ by
	\begin{equation*}
	\tilde{D}_{k,r}{(x)}=\frac{\cos \left(k+\frac{r}{2}\right)x}{2\sin \frac{rx}{2}}, \hspace{15mm} 
	D_{k,r}{(x)}=\frac{\sin \left(k+\frac{r}{2}\right)x}{2\sin \frac{rx}{2}}
	\end{equation*}
	the Dirichlet type kernels.

	\begin{lem} \text{(\cite{Szal}, \cite{Szal1}} \label{roznicowanie})
		Let $r,m,n \in \mathbb{N}, l\in \mathbb{Z}$ and $(a_k) \subset \mathbb{C}$. If $x \neq \frac{2l\pi}{r}$, then for $m \geq n$
		\begin{equation}\label{sinus_roznicowanie}
		\sum\limits_{k=n}^{m} a_{k} sin(kx) = 
		- \sum\limits_{k=n}^{m} \Delta_{r} a_{k} \tilde{D}_{k,r}(x)
		+
		\sum\limits_{k=m+1}^{m+r} a_{k} \tilde{D}_{k,-r}(x) 
		+
		\sum\limits_{k=n}^{n+r-1} a_{k} \tilde{D}_{k,-r}(x) 
		\end{equation}
		and
		\begin{equation}\label{cosinus_roznicowanie}
		\sum\limits_{k=n}^{m}a_{k}\cos kx = \sum\limits_{k=n}^{m} \Delta_r a_k \tilde{D}_{k,r}(x) - \sum\limits_{k=m+1}^{m+r} a_k \tilde{D}_{k,-r}(x) + \sum\limits_{k=n}^{n+r-1} a_k \tilde{D}_{k,-r}(x)
		\end{equation}
	\end{lem}

	\begin{lem}\text{(\cite{Kubiak})} \label{roznicowanie_exp}
		Niech $r,m,n\in \mathbb{N}, l \in \mathbb{Z}$ oraz $a=\{c_k\}_ \subset \mathbb{C}$. If $x \neq \frac{2l\pi}{r}$, then for  $m \geq n$
		\begin{equation*}
		\sum\limits_{k=n}^{m}a_{k}e^{ikx} = \frac{-i}{2\sin\left(\frac{rx}{2}\right)}
		\left( \sum\limits_{k=n}^{m} \Delta_r a_k e^{-i\left( k+\frac{r}{2} \right)x} 
		- \sum\limits_{k=m+1}^{m+r}  e^{-i\left( k-\frac{r}{2} \right)x}
		+ \sum\limits_{k=n}^{n+r-1}  e^{-i\left( k-\frac{r}{2} \right)x}
		\right).
		\end{equation*}
	\end{lem}

	\begin{lem}\label{lemat_do_warunku_z_logarytmem}
		Let $n, N \in \mathbb{N}$. Then for $p \geq 1$
		\begin{equation*}
		\int_{n+N^{\frac{1}{p}}}^{n+N} \frac{1}{k\ln k} dk \leq \ln p.
		\end{equation*}
		
		\begin{proof}
			This inequality is true for $p =1 $. Consider the function $f(p) = \left(n+N^{\frac{1}{p}}\right)^p$ for $p > 0$. We get:
			\begin{equation*}
			f'(p) = p\left(n+N^{\frac{1}{p}}\right)^{p-1} \frac{1}{p}N^{\frac{1}{p}-1} =  \left(n+N^{\frac{1}{p}}\right)^{p-1} N^{\frac{1}{p}-1} \geq 0 \text{ for all }p>0.
			\end{equation*}
			It means that the function is non-decresing with respect to $p$. 
			\newline Thus:
			\begin{equation}
			n+N  = f(1) \leq f(p) = \left(n+N^{\frac{1}{p}}\right)^p \text{ for } p\geq 1.
			\end{equation}
			Hence we get that:
			\begin{equation}\label{nierowosclemat_z_calka}
			\ln\left(n+N\right) \leq \ln\left(n+N^{\frac{1}{p}}\right)^p.
			\end{equation}
			Therefore, integrating by substitution with $\ln k = t$ and using (\ref{nierowosclemat_z_calka}), we get
			\begin{equation*}
			\int_{n+N^{\frac{1}{p}}}^{n+N} \frac{1}{k \ln k } dk =
		    \int_{\ln(n+N)^{\frac{1}{p}}}^{\ln(n+N)} \frac{1}{t} dt
		    = \ln(\ln(n+N)) - \ln(\ln(n+N^\frac{1}{p}))
			\end{equation*}
			\begin{equation*}
			 = \ln\left(\frac{\ln(n+N)}{\ln\left(n+N^{\frac{1}{p}}\right)}\right)
			= \ln\left(p \frac{\ln\left(n+N\right)}{\ln\left(n+N^{\frac{1}{p}}\right)^p}\right)
			\leq \ln p
			\end{equation*}
			and the proof is completed.
		\end{proof}
	\end{lem}

\section{Proofs od the main results}

\subsection{Proof of the theorem 10}
	Let $\epsilon >0$. Then from \eqref{warunek_n_ln_n} and \eqref{zbiezne_w_punkach} we obtain: 
	\begin{equation}\label{prof12_2}
	n\ln n \left\vert b_n \right\vert <\varepsilon,
	\end{equation}
	\begin{equation}\label{4.1}
	\left\vert \sum\limits_{k=n}^{\infty } b_{k}\sin(k\frac{2l\pi}{r})\right\vert < \varepsilon,
	\end{equation}%
	and
	\begin{equation}\label{4.2}
	\left\vert \sum\limits_{k=n}^{n+N } b_{k}\sin(k\frac{2l\pi}{r})\right\vert < \varepsilon,
	\end{equation}
	for all $n > N_\epsilon$ and $N \in \mathbb{N}$, where $l= 1,...,[\frac{r}{2}]-1 $ when $r$ is an even number and $l= 1,...,[\frac{r}{2}] $ when $r$ is an odd number.
	Denote by
	\begin{equation*}
	\tau_{n}\left( x\right) =\sum\limits_{k=n}^{\infty }b_{k}\sin (ckx).
	\end{equation*}%
	We will show that 
	\begin{equation}\label{4.3}
	\left\vert \tau_{n}\left( x\right) \right\vert \ll \varepsilon  
	\end{equation}%
	holds for any $n\geq \max\{N_\varepsilon,2\}$ and $x \in \mathbb{R}$. Since $\tau_{n}\left( 0 \right) =0$ and $\tau_{n}\left( \frac{\pi}{c} \right) =0$ it suffices to prove \eqref{4.3} for $0<x< \frac{\pi}{c} $.
	
	First, we will show that \eqref{4.3} is valid for $x=\frac{2l\pi }{rc}$, where $l$
	is an integer number such that $0<2l<r.$ Using \eqref{4.1} we get
	\begin{equation}\label{4.0}
	\left\vert \tau_{n}\left( \frac{2l\pi }{rc}\right) \right\vert  < \varepsilon.
	\end{equation}%
	
	Now, we prove that \eqref{4.3} holds for $\frac{2l\pi }{rc}<x\leq \frac{2l\pi 
	}{rc}+\frac{\pi }{rc}$, where $0\leq 2l<r$.
	
	Let $N^{\frac{1}{p}}:=N^{\frac{1}{p}}\left( x\right) \geq r $ be a natural number such that%
	\begin{equation}
	\frac{2l\pi }{rc}+\frac{\pi }{c(N+1)^{\frac{1}{p}}}<x\leq \frac{2l\pi }{rc}+\frac{\pi }{cN^{\frac{1}{p}}}.
	\label{p4}
	\end{equation}%
	Then%
	\begin{equation*}
	\tau_{n}\left( x\right) =
	\sum\limits_{k=n}^{n+N^{\frac{1}{p}}-1}b_{k}\sin(ckx)
	+ \sum\limits_{k=n+N^{\frac{1}{p}}}^{n+N}b_{k}\sin (ckx)
	+ \sum\limits_{k=n+N+1}^{\infty }b_{k}\sin (ckx)
	\end{equation*}
	\begin{equation*}
	=\tau_{n}^{\left( 1\right) }\left(
	x\right)  +\tau_{n}^{\left( 2\right) }\left( x\right) +\tau_{n}^{\left( 3\right) }\left( x\right) .
	\end{equation*}%
	Applying Lagrange's mean value theorem to the function $f\left( x\right)
	=\sin (ckx)$ on the interval $\left[ \frac{2l\pi }{rc},x\right] $ we obtain that for each $k$
	there exists $y_k\in \left( \frac{2l\pi }{rc},x\right) $ such that  %
	\begin{equation*}
	\sin (ckx)-\sin \left( k\frac{2l\pi }{r}\right) =ck\cos (cky_k)\left( x-\frac{2l\pi }{%
		rc}\right) .
	\end{equation*}%
	Using this we get%
	\begin{equation*}
	\tau_{n}^{\left( 1\right) }\left( x\right)
	=\sum\limits_{k=n}^{n+N^{\frac{1}{p}}-1}ckb_{k}\cos (cky_k)\left( x-\frac{2l\pi }{rc}\right)
	+\sum\limits_{k=n}^{n+N^{\frac{1}{p}}-1}b_{k}\sin \left( k\frac{2l\pi }{r}\right)
	\end{equation*}%
	\begin{equation*}
	=  \tau_n^{(1.1)}(x) + \tau_n^{(1.2)}.
	\end{equation*}%
	From \eqref{4.2} we have
	\begin{equation*}
	\left\vert \tau_n^{(1.2)} \right\vert <\varepsilon.
	\end{equation*}
	
	By \eqref{p4} and \eqref{prof12_2}%
	\begin{equation*}
	\left\vert \tau_{n}^{\left( 1.1\right) }\left( x\right) \right\vert \leq \left(
	x-\frac{2l\pi }{rc}\right) \sum\limits_{k=n}^{n+N^{\frac{1}{p}}-1}ck\left\vert b_{k}\right\vert 
	\leq \left(x-\frac{2l\pi}{rc}\right) \sum\limits_{k=n}^{n+N^{\frac{1}{p}}-1}\frac{ck\ln k}{\ln k}\left\vert b_{k}\right\vert 
	\end{equation*}
	\begin{equation}
	< \left(x-\frac{2l\pi}{rc}\right) \sum\limits_{k=n}^{n+N^{\frac{1}{p}}-1}\frac{c\varepsilon}{\ln k}
	\leq \frac{\pi \varepsilon}{\ln 2} .  \label{p10}
	\end{equation}%
	
	Using Lemma \ref{lemat_do_warunku_z_logarytmem} we have
	\begin{equation*}
	\left\vert \tau_n^{(2)}(x)\right\vert = \left\vert \sum\limits_{k=n+N^{\frac{1}{p}}}^{n+N}b_{k}\sin (ckx) \right\vert
	\leq \sum\limits_{k=n+N^{\frac{1}{p}}}^{n+N}\frac{k \ln k }{k \ln k} \left\vert b_k \right\vert
	\ll  \varepsilon \int_{n+N^{\frac{1}{p}}}^{n+N} \frac{1}{k \ln k} dk 
	\leq  \varepsilon \ln p
	\end{equation*}
	
	If $\left( b_{n}\right) \in GM\left(p,{}_3\beta(q),r\right) $, then using Lemma \ref{roznicowanie}, we get
	\begin{equation*}
	\left\vert \tau_{n}^{\left( 3\right) }\left( x\right) \right\vert =\left\vert
	\sum\limits_{j=0}^{\infty }\sum\limits_{k=2^{j}\left( n+N+1\right)
	}^{2^{j+1}\left( n+N+1\right) -1}b_{k}\sin (ckx)\right\vert
	\end{equation*}%
	\begin{equation*}
	\leq \sum\limits_{j=0}^{\infty }\left\vert \frac{-1}{2\sin \left(
		crx/2\right) }\left\{\sum\limits_{k=2^{j}\left( n+N+1\right) }^{2^{j+1}\left(
		n+N+1\right) -1}\left( b_{k}-b_{k+r}\right) \cos \left( k+\frac{r}{2}\right)c
	x\right.\right.
	\end{equation*}%
	\begin{equation*}
	\left. \left. +\sum\limits_{k=2^{j+1}\left( n+N+1\right) }^{2^{j+1}\left(
		n+N+1\right) +r-1}b_{k}\cos \left( k-\frac{r}{2}\right)c
	x-\sum\limits_{k=2^{j}\left( n+N+1\right) }^{2^{j}\left( n+N+1\right)
		+r-1}b_{k}\cos \left( k-\frac{r}{2}\right) cx\right\} \right\vert
	\end{equation*}%
	\begin{equation*}
	\leq \frac{1}{2\left\vert \sin \left( crx/2\right) \right\vert }%
	\sum\limits_{j=0}^{\infty }\left\{ \sum\limits_{k=2^{j}\left( n+N+1\right)
	}^{2^{j+1}\left( n+N+1\right) -1}\left\vert b_{k}-b_{k+r}\right\vert
	+\sum\limits_{k=2^{j+1}\left( n+N+1\right) }^{2^{j+1}\left( n+N+1\right)
		+r-1}\left\vert b_{k}\right\vert \right. 
	\end{equation*}%
	\begin{equation*}
	\left. +\sum\limits_{k=2^{j}\left( n+N+1\right)
	}^{2^{j}\left( n+N+1\right) +r-1}\left\vert b_{k}\right\vert \right\}.
	\end{equation*}
	Further applaying the H\"{o}lder inequality with $p>1$, the inequality 
	$\frac{rc}{\pi }x-2l\leq
	\left\vert \sin \frac{rcx}{2}\right\vert $ $\left( x\in \left[ \frac{2l\pi }{rc%
	},\frac{2l\pi }{rc}+\frac{\pi }{rc}\right] \text{ and }0\leq 2l<r\right)$, \eqref{p4} and \eqref{prof12_2}, we obtain
	\begin{equation*}
	\left\vert \tau_{n}^{\left( 3\right) }\left( x\right) \right\vert \leq \frac{1}{\frac{rc}{\pi }x-2l}\sum\limits_{j=0}^{\infty }\left\{ \left(
	\sum\limits_{k=2^{j}\left( n+N+1\right) }^{2^{j+1}\left( n+N+1\right)
		-1}\left\vert b_{k}-b_{k+r}\right\vert^p\right)^{\frac{1}{p}} \left( \sum\limits_{k=2^{j}\left( n+N+1\right) }^{2^{j+1}\left( n+N+1\right)
		-1}1\right) ^{1-\frac{1}{p}}\right.
	\end{equation*}
	\begin{equation*}
	\left.+\sum\limits_{k=2^{j}\left(
		n+N+1\right) }^{2^{j}\left( n+N+1\right) +r-1}\left\vert b_{k}\right\vert
	+\sum\limits_{k=2^{j+1}\left(
		n+N+1\right) }^{2^{j+1}\left( n+N+1\right) +r-1}\left\vert b_{k}\right\vert
	\right\}
	\end{equation*}%
	\begin{equation*}
	\leq \frac{(N+1)^\frac{1}{p}}{r}\sum\limits_{j=0}^{\infty }\left\{
	\frac{C\left( 2^j(n+N+1)\right)^{1-\frac{1}{p}}}{2^j(n+N+1)} \underset{m\geq b\left(2^j(n+N+1)\right)}{\sup }m \left( \frac{1}{m} \sum\limits_{k=m}^{2m-1 } \left\vert b_k\right\vert^q\right)^{\frac{1}{q}} 
	\right.
	\end{equation*}
	
	\begin{equation*}
	\left.
	+\sum\limits_{k=2^{j}\left(
		n+N+1\right) }^{2^{j}\left( n+N+1\right) +r-1}\left\vert b_{k}\right\vert
	+\sum\limits_{k=2^{j+1}\left(
		n+N+1\right) }^{2^{j+1}\left( n+N+1\right) +r-1}\left\vert b_{k}\right\vert
	\right\}
	\end{equation*}%
	\begin{equation*}
	= \frac{(N+1)^\frac{1}{p}}{r}\sum\limits_{j=0}^{\infty }\left\{
	\frac{C}{\left(2^j(n+N+1)\right)^{\frac{1}{p}}} \underset{m\geq b\left(2^j(n+N+1)\right)}{\sup }m^{1-\frac{1}{q}} \left( \sum\limits_{k=m}^{2m-1 } \left( \frac{k \ln k\left\vert b_k \right\vert}{k\ln k}\right)^q\right)^{\frac{1}{q}}
	\right.
	\end{equation*}
	\begin{equation*}
	\left.
	+\sum\limits_{k=2^{j}\left(
		n+N+1\right) }^{2^{j}\left( n+N+1\right) +r-1}\frac{k\ln k\left\vert b_k \right\vert}{k\ln k}
	+\sum\limits_{k=2^{j+1}\left(
		n+N+1\right) }^{2^{j+1}\left( n+N+1\right) +r-1}\frac{k\ln k\left\vert b_k \right\vert}{k\ln k}
	\right\}
	\end{equation*}
	\begin{equation*}
	< \frac{\varepsilon (N+1)^\frac{1}{p}}{r}\sum\limits_{j=0}^{\infty }\left\{
	\frac{C}{\left(2^j(n+N+1)\right)^{\frac{1}{p}}} \underset{m\geq b\left(2^j(n+N+1)\right)}{\sup }m^{1-\frac{1}{q}} \left( \sum\limits_{k=m}^{2m-1 } \left( \frac{1 }{k}\right)^q\right)^{\frac{1}{q}}
	\right.
	\end{equation*}
	\begin{equation*}
	\left.
	+\sum\limits_{k=2^{j}\left(
		n+N+1\right) }^{2^{j}\left( n+N+1\right) +r-1}\frac{1}{k}
	+\sum\limits_{k=2^{j+1}\left(
		n+N+1\right) }^{2^{j+1}\left( n+N+1\right) +r-1}\frac{1}{k}
	\right\}
	\end{equation*}
	\begin{equation*}
	\leq \frac{\varepsilon (N+1)^\frac{1}{p}}{r}\sum\limits_{j=0}^{\infty }\left\{
	\frac{C}{\left(2^j(n+N+1)\right)^{\frac{1}{p}}} \underset{m\geq b\left(2^j(n+N+1)\right)}{\sup } \left(m^{1-\frac{1}{q}} m^{-1} m^{\frac{1}{q}}\right)
	\right.
	\end{equation*}
	\begin{equation*}
	\left.
	+ \frac{3}{2}r\left( 2^j(n+N+1)\right)^{-\frac{1}{p}}
	\right\}.
	\end{equation*}
	Elementary calculations give:
	\begin{equation*}
	\left\vert \tau_{n}^{\left( 3\right) }\left( x\right) \right\vert < \frac{\varepsilon (N+1)^\frac{1}{p}}{r}\sum\limits_{j=0}^{\infty }\left\{
	\left(2^j(n+N+1)\right)^{-\frac{1}{p}}  \left( C+\frac{3}{2}r\right) \right\}
	\leq \frac{\varepsilon (N+1)^\frac{1}{p}\left( C+\frac{3}{2}r\right)}{r \left((n+N+1)\right)^{\frac{1}{p}}}\sum\limits_{j=0}^{\infty }\left( \frac{1}{2^\frac{1}{p}}  \right)^j
	\end{equation*}
	\begin{equation}
	\leq \frac{\varepsilon \left( C+\frac{3}{2}r\right)}{r} \frac{1}{1-2^{\frac{1}{p}}}.  \label{4.6}
	\end{equation}
	
	Finally, we prove that \eqref{4.3} is true for $\frac{2l\pi }{rc}+\frac{\pi }{rc%
	}\leq x<\frac{2\left( l+1\right) \pi }{rc}$, where $0<2\left( l+1\right) \leq
	r.$
	
	Let $M^{\frac{1}{p}}:=M^{\frac{1}{p}}\left( x\right) \geq r$ be a natural number such that%
	\begin{equation}
	\frac{2\left( l+1\right) \pi }{rc}-\frac{\pi }{cM^{\frac{1}{p}}}\leq x<\frac{2\left(
		l+1\right) \pi }{rc}-\frac{\pi }{c(M+1)^{\frac{1}{p}}}.  \label{p8}
	\end{equation}%
	Then%
	\begin{equation*}
	\tau_{n}\left( x\right) =
	\sum\limits_{k=n}^{n+M^{\frac{1}{p}}-1}b_{k}\sin(ckx)
	+ \sum\limits_{k=n+M^{\frac{1}{p}}}^{n+M}b_{k}\sin (ckx)
	+ \sum\limits_{k=n+M+1}^{\infty }b_{k}\sin (ckx)
	\end{equation*}
	\begin{equation*}
	=\tau_{n}^{\left( 4\right) }\left(
	x\right)  +\tau_{n}^{\left( 5\right) }\left( x\right) +\tau_{n}^{\left( 6\right) }\left( x\right) .
	\end{equation*}%
	Applying Lagrange's mean value theorem to the function $f\left( x\right)
	=\sin (ckx)$ on the interval $\left[ x,\frac{2\left( l+1\right) \pi }{rc}\right] 
	$ we obtain that for each $k$ there exists $z_k\in \left( x,\frac{2\left( l+1\right) \pi }{rc%
	}\right) $ such that%
	\begin{equation*}
	\sin \left( k\frac{2\left( l+1\right) \pi }{r}\right) -\sin (ckx)=ck\cos
	(ckz_k)\left( \frac{2\left( l+1\right) \pi }{rc}-x\right) .
	\end{equation*}%
	Using this we get%
	\begin{equation*}
	\tau_{n}^{\left( 4\right) }\left( x\right)
	=\sum\limits_{k=n}^{n+M^{\frac{1}{p}}-1}ckb_{k}\cos (ckz_k)\left( \frac{2\left( l+1\right) \pi }{rc}-x\right) +\sum\limits_{k=n}^{n+M^{\frac{1}{p}}-1}b_{k}\sin \left( k\frac{2\left(
		l+1\right) \pi }{r}\right)
	\end{equation*}%
	\begin{equation*}
	=\tau_{n}^{\left( 4.1\right) }\left( x\right) +\tau_{n}^{\left( 4.2\right) }.
	\end{equation*}%
	From \eqref{4.1} we have 
	\begin{equation*}
	\left\vert \tau_n^{(4.2)} \right\vert <\varepsilon.
	\end{equation*}
	By \eqref{p8} and \eqref{prof12_2}%
	\begin{equation}
	\left\vert \tau_{n}^{\left( 4.1\right) }\left( x\right) \right\vert \leq \left( 
	\frac{2\left( l+1\right) \pi }{rc}-x\right)
	\sum\limits_{k=n}^{n+M^{\frac{1}{p}}-1}ck\left\vert b_{k}\right\vert \leq \frac{\pi}{M^{\frac{1}{p}}} \sum\limits_{k=n}^{n+M^{\frac{1}{p}}-1}\frac{k\ln k}{\ln k}\left\vert b_{k}\right\vert 
	\leq \frac{\pi \varepsilon}{\ln 2}  \label{p14}
	\end{equation}%

	If $\left( b_{n}\right) \in GM\left(p,{}_3 \beta(q) ,r\right) $, we obtain%
	\begin{equation*}
	\left\vert \tau_{n}^{\left( 6\right) }\left( x\right) \right\vert =\left\vert
	\sum\limits_{j=0}^{\infty }\sum\limits_{k=2^{j}\left( n+M+1\right)
	}^{2^{j+1}\left( n+M+1\right) -1}b_{k}\sin (ckx)\right\vert
	\end{equation*}%
	\begin{equation*}
	\leq \sum\limits_{j=0}^{\infty }\left\vert \frac{-1}{2\sin \left(
		crx/2\right) }\left\{ \sum\limits_{k=2^{j}\left( n+M+1\right) }^{2^{j+1}\left(
		n+M+1\right) -1}\left( b_{k}-b_{k+r}\right) \cos \left( \left( k+\frac{r}{2}\right)c
	x\right)\right.\right.
	\end{equation*}%
	\begin{equation*}
	\left. \left. +\sum\limits_{k=2^{j+1}\left( n+M+1\right) }^{2^{j+1}\left(
		n+M+1\right) +r-1}b_{k}\cos \left( k-\frac{r}{2}\right)c
	x-\sum\limits_{k=2^{j}\left( n+M+1\right) }^{2^{j}\left( n+M+1\right)
		+r-1}b_{k}\cos \left(\left( k-\frac{r}{2}\right) cx\right)\right\} \right\vert.
	\end{equation*}%
	\begin{equation*}
	\leq \frac{1}{2\left\vert \sin \left( crx/2\right) \right\vert }%
	\sum\limits_{j=0}^{\infty }\left\{ \sum\limits_{k=2^{j}\left( n+M+1\right)
	}^{2^{j+1}\left( n+M+1\right) -1}\left\vert b_{k}-b_{k+r}\right\vert
	+\sum\limits_{k=2^{j+1}\left( n+M+1\right) }^{2^{j+1}\left( n+M+1\right)
		+r-1}\left\vert b_{k}\right\vert \right.
	\end{equation*}%
	\begin{equation*}
	\left. +\sum\limits_{k=2^{j}\left( n+M+1\right)}^{2^{j}\left( n+M+1\right) +r-1}\left\vert b_{k}\right\vert \right\}.
	\end{equation*}
	Next, applying  the H\"{o}lder inequality with $p>1$, then using Lemma \ref{roznicowanie}, the inequality 
	\newline $2\left( l+1\right) -\frac{rc}{%
		\pi }x\leq \left\vert \sin \frac{rcx}{2}\right\vert $ $\left( x\in \left[ 
	\frac{2l\pi }{rc}+\frac{\pi }{rc},\frac{2\left( l+1\right) \pi }{rc}\right] 
	\text{ and }0<2\left( l+1\right) \leq r\right) $, \eqref{p8} and \eqref{prof12_2}, we get
	\begin{equation*}
	\left\vert \tau_{n}^{\left( 6\right) }\left( x\right) \right\vert \leq \frac{1}{2(l+1)-\frac{rc}{\pi}x} \sum\limits_{j=0}^{\infty }\left\{ \left(
	\sum\limits_{k=2^{j}\left( n+M+1\right) }^{2^{j+1}\left( n+M+1\right)
		-1}\left\vert b_{k}-b_{k+r}\right\vert^p\right)^{\frac{1}{p}}  \left(\sum\limits_{k=2^{j}\left( n+M+1\right) }^{2^{j+1}\left( n+M+1\right)
		-1} 1\right)^{1-\frac{1}{p}}
	\right.
	\end{equation*}
	\begin{equation*}
	\left.
	+\sum\limits_{k=2^{j}\left(
		n+M+1\right) }^{2^{j}\left( n+M+1\right) +r-1}\left\vert b_{k}\right\vert
	+\sum\limits_{k=2^{j+1}\left(
		n+M+1\right) }^{2^{j+1}\left( n+M+1\right) +r-1}\left\vert b_{k}\right\vert
	\right\}
	\end{equation*}%
	\begin{equation*}
	\leq \frac{(M+1)^\frac{1}{p}}{r}\sum\limits_{j=0}^{\infty }\left\{
	\frac{C\left( 2^j(n+M+1)\right)^{1-\frac{1}{p}}}{2^j(n+M+1)} \underset{m\geq b\left(2^j(n+M+1)\right)}{\sup }m \left( \frac{1}{m} \sum\limits_{k=m}^{2m-1 } \left\vert b_k\right\vert^q \right)^{\frac{1}{q}} 
	\right.
	\end{equation*}
	\begin{equation*}
	\left.
	+\sum\limits_{k=2^{j}\left(
		n+M+1\right) }^{2^{j}\left( n+M+1\right) +r-1}\left\vert b_{k}\right\vert
	+\sum\limits_{k=2^{j+1}\left(
		n+M+1\right) }^{2^{j+1}\left( n+M+1\right) +r-1}\left\vert b_{k}\right\vert
	\right\}
	\end{equation*}%
	\begin{equation*}
	= \frac{(M+1)^\frac{1}{p}}{r}\sum\limits_{j=0}^{\infty }\left\{
	\frac{C}{\left(2^j(n+M+1)\right)^{\frac{1}{p}}} \underset{m\geq b\left(2^j(n+M+1)\right)}{\sup }m^{1-\frac{1}{q}} \left( \sum\limits_{k=m}^{2m-1 } \left( \frac{k\ln k\left\vert b_k \right\vert}{k\ln k}\right)^q\right)^{\frac{1}{q}} \right.
	\end{equation*}
	\begin{equation*}
	\left. +\sum\limits_{k=2^{j}\left(
		n+M+1\right) }^{2^{j}\left( n+M+1\right) +r-1}\frac{k\ln k\left\vert b_k \right\vert}{k\ln k}
	+\sum\limits_{k=2^{j+1}\left(
		n+M+1\right) }^{2^{j+1}\left( n+M+1\right) +r-1}\frac{k\ln k\left\vert b_k \right\vert}{k\ln k}
	\right\}.
	\end{equation*}%
	\begin{equation*}
	\end{equation*}
	\begin{equation*}
	< \frac{\varepsilon (M+1)^\frac{1}{p}}{r}\sum\limits_{j=0}^{\infty }\left\{
	\frac{C}{\left(2^j(n+M+1)\right)^{\frac{1}{p}}} \underset{m\geq b\left(2^j(n+M+1)\right)}{\sup }m^{1-\frac{1}{q}} \left( \sum\limits_{k=m}^{2m-1 } \left( \frac{1 }{k}\right)^q\right)^{\frac{1}{q}}
	\right.
	\end{equation*}
	\begin{equation*}
	\left.
	+\sum\limits_{k=2^{j}\left(
		n+M+1\right) }^{2^{j}\left( n+M+1\right) +r-1}\frac{1}{k}
	+\sum\limits_{k=2^{j+1}\left(
		n+M+1\right) }^{2^{j+1}\left( n+M+1\right) +r-1}\frac{1}{k}
	\right\}.
	\end{equation*}
	Elementary calculations give
	\begin{equation*}
	\left\vert \tau_{n}^{\left( 6\right) }\left( x\right) \right\vert < \frac{\varepsilon (M+1)^\frac{1}{p}}{r}\sum\limits_{j=0}^{\infty }\left\{
	\frac{C}{\left(2^j(n+M+1)\right)^{\frac{1}{p}}} \underset{m\geq b\left(2^j(n+M+1)\right)}{\sup }\left(m^{1-\frac{1}{q}} m^{-1} m^{\frac{1}{q}}\right)
	\right.
	\end{equation*}
	\begin{equation*}
	\left.
	+ r \left( 2^j(n+M+1)\right)^{-\frac{1}{p}}
	+ \frac{1}{2}r\left( 2^j(n+M+1)\right)^{-\frac{1}{p}}
	\right\}
	\end{equation*}
	\begin{equation*}
	\leq \frac{\varepsilon (M+1)^\frac{1}{p}}{r}\sum\limits_{j=0}^{\infty }\left\{
	\left(2^j(n+M+1)\right)^{-\frac{1}{p}}  \left( C+\frac{3}{2}r\right) \right\}
	\end{equation*}
	\begin{equation}
	\leq \frac{\varepsilon (M+1)^\frac{1}{p}\left( C+\frac{3}{2}r\right)}{r \left((n+M+1)\right)^{\frac{1}{p}}}\sum\limits_{j=0}^{\infty }\left( \frac{1}{2^\frac{1}{p}}  \right)^j
	\leq \frac{\varepsilon \left( C+\frac{3}{2}r\right)}{r} \frac{1}{1-2^{\frac{1}{p}}}.  \label{4.6}
	\end{equation}
	
	Joining the obtained estimates the uniform convergence of
	series \eqref{sinus} follows and thus the proof is complete. $\square $
	
\subsection{Proof of Remark 1}
	Let $p>1$, $c>1$ and for $n \in \mathbb{N}$:
	\begin{center}
		$a_n =\left\{ \begin{array}{l} \frac{3}{n\ln(n+1)},\text{ when } n=1\text{ (mod }3),
		\\\frac{1}{n\ln(n+1)},\text{  when } n=2\text{ (mod }3),
		\\\frac{1}{n\ln(n+1)},\text{  when } n=0 (\mod 3) \text{ and } n\neq 0\text{ (mod }6),
		\\\frac{1}{(n-3)\ln(n-2)} + \frac{1}{n^{1+\frac{1}{p}}\ln(n+1)},\text{  when } n=0\text{ (mod }6). \end{array} \right.$
	\end{center}
	First, we prove that $(a_n) \in GM(p,{}_3\beta(1),3)$. Let
	\begin{equation*}
	A_n = \{ k \in \mathbb{N} : n\leq k \leq 2n-1\text{ and } k =1 (\mod 3) \},
	\end{equation*}
	\begin{equation*}
	B_n = \{ k \in \mathbb{N} : n\leq k\leq 2n-1 \text{ and } k=1 (\mod 3) \},
	\end{equation*}
	\begin{equation*}
	C_n = \{ k \in \mathbb{N} : n\leq k\leq 2n-1\text{ and } k=0 (\mod 3) \text{ and } k\neq 0 (\mod 6) \}.
	\end{equation*}
	\begin{equation*}
	D_n = \{ k \in \mathbb{N} : n\leq k\leq 2n-1\text{ and } k=0 (\mod 6) \}.
	\end{equation*}
	We have:
	\begin{equation*}
	\left\{ \sum_{k=n}^{2n-1} \left| a_k - a_{k+3} \right|^p \right\}^{\frac{1}{p}}
	= \left\{ \sum_{k\in A_n} \left| a_k - a_{k+3} \right|^p
	+ \sum_{k\in B_n} \left| a_k - a_{k+3} \right|^p
	+ \sum_{k\in C_n} \left| a_k - a_{k+3} \right|^p
	\right.
	\end{equation*}
	\begin{equation*}
	\left.
	+ \sum_{k\in D_n} \left| a_k - a_{k+3} \right|^p \right\}^{\frac{1}{p}}
	= \left\{
	\sum_{k\in A_n} \left| \frac{3}{k\ln(k+1)} - \frac{3}{(k+3)\ln(k+4)} \right|^p 
	\right.
	\end{equation*}
	\begin{equation*}
	+ \sum_{k\in B_n} \left| \frac{1}{k\ln(k+1)} - \frac{1}{(k+3)\ln(k+1)} \right|^p
	\end{equation*}
	\begin{equation*}
	+ \sum_{k\in C_n} \left| \frac{1}{k\ln(k+1)} - \frac{1}{k\ln(k+1)} - \frac{1}{(k+3)^{1+\frac{1}{p}}\ln(k+4)} \right|^p
	\end{equation*}
	\begin{equation*}
	\left.
	+ \sum_{k\in D_n} \left| \frac{1}{(k-3)\ln(k-2)} + \frac{1}{k^{1+\frac{1}{p}}\ln(k+1)} - \frac{1}{(k+3)\ln(k+4)} \right|^p
	\right\}^{\frac{1}{p}}
	\end{equation*}
	\begin{equation*}
	= \left\{
	\sum_{k\in A_n} 3^p \left| \frac{1}{k\ln(k+1)} - \frac{1}{(k+3)\ln(k+4)} \right|^p
	+  \sum_{k\in B_n} \left| \frac{1}{k\ln(k+1)} - \frac{1}{(k+3)\ln(k+1)} \right|^p
	\right.
	\end{equation*}
	\begin{equation*}
	+ \sum_{k\in C_n} \left|  \frac{1}{(k+3)^{1+\frac{1}{p}}\ln(k+4)} \right|^p
	\end{equation*}
	\begin{equation*}
	\left.
	+ \sum_{k\in D_n} \left(\left| \frac{1}{(k-3)\ln(k-2)}  - \frac{1}{(k+3)\ln(k+4)} \right| + \frac{1}{k^{1+\frac{1}{p}}\ln(k+1)}\right)^p
	\right\}^{\frac{1}{p}}.
	\end{equation*}
	\begin{equation*}
	= \left\{
	\sum_{k\in A_n} 3^p S_1^p
	+ \sum_{k\in B_n} S_1
	+ \sum_{k\in C_n} \frac{1}{(k+3)^{1+\frac{1}{p}}\ln(k+4)} 
	+ \sum_{k\in D_n} \left( S_2 + \frac{1}{k^{1+\frac{1}{p}}\ln(k+1)} \right)^p
	\right\}^{\frac{1}{p}}
	\end{equation*}
	For $S_1$, we have
	\begin{equation*}
	S_1 = \left| \frac{1}{k\ln(k+1)} - \frac{1}{(k+3)\ln(k+4)} \right| =  \frac{\left| (k+3)\ln(k+4)-k\ln(k+1) \right|}{k(k+3)\ln(k+1)\ln(k+4)} 
	\end{equation*}
	\begin{equation*}
	\leq \frac{\frac{3k}{k+1}+3\ln(k+3)}{k(k+3)\ln(k+1)\ln(k+4)}
	\leq \frac{6\ln(k+4)}{k^2\ln(k+1)\ln(k+4)} = \frac{6}{k^2\ln(k+1)}
	\end{equation*}
	Further, we obtain
	\begin{equation*}
	S_2 = \left| \frac{1}{(k-3)\ln(k-2)}  - \frac{1}{(k+3)\ln(k+4)} \right| 
	= \frac{\left|(k+3)\ln(k+4)-(k-3)\ln(k-2)\right|}{(k-3)(k+3)\ln(k-2)\ln(k+4)}
	\end{equation*}
	\begin{equation*}
	\leq
	\frac{(k-3)\left| \ln(k+4)-\ln(k-2)\right| + 6\ln(k+4)}{(k-3)(k+3)\ln(k-2)\ln(k+4)}
	\end{equation*}
	Using elementary calculations and Lagrange theorem, we have
	\begin{equation*}
	\left| \ln(k+4)-\ln(k-2) \right| \leq \frac{1}{k-2} \left( k+4-k+2 \right) = \frac{6}{k-2}
	\end{equation*}
	So, we obtain
	\begin{equation*}
	S_2 \leq \frac{\frac{6(k-3)}{k-2}+ 6\ln(k+4)}{(k-3)(k+3)\ln(k-2)\ln(k+4)}
	\leq  \frac{12}{(k-3)(k+3)\ln(k-2)}
	\leq \frac{48}{k^2\ln(k+1)}
	\end{equation*}
	Hence, we have
	\begin{equation*}
	\left\{ \sum_{k=n}^{2n-1} \left| a_k - a_{k+3} \right|^p \right\}^{\frac{1}{p}}
	\leq \left\{ \sum_{k\in A_n} 3^p \left( \frac{6}{k^2\ln(k+1)} \right)^p
	+ \sum_{k\in B_n} \left( \frac{6}{k^2\ln(k+1)} \right)^p
	\right.
	\end{equation*}
	\begin{equation*}
	\left.
	+ \sum_{k\in C_n} \left( \frac{1}{k^{1+\frac{1}{p}}\ln(k+1)} \right)^p 
	+ \sum_{k\in D_n} \left( \frac{48}{k^2\ln(k+1)} + \frac{1}{k^{1+\frac{1}{p}}\ln(k+1)} \right)^p \right\}^\frac{1}{p}
	\end{equation*}
	\begin{equation*}
	\leq 49 \left\{ \sum_{k=n}^{2n-1} \left( \frac{1}{k^{1+\frac{1}{p}}\ln(k+1)} \right)^p \right\}^{\frac{1}{p}}
	\leq 49 \frac{1}{n^{1+\frac{1}{p}}\ln(n+1)} n^{\frac{1}{p}} = 49 \frac{1}{n\ln(n+1)}
	\end{equation*}
	\begin{equation*}
	\leq 147 \frac{1}{n} \sum_{k=n}^{2n} |a_n|
	\leq 147 {}_3\beta(1)
	\end{equation*}
	Hence $(a_n) \in GM(p,{}_3\beta(1),3)$. Now, we will show that $(a_n) \notin GM(1,{}_3\beta(1),3)$.
	\begin{equation*}
	\sum_{k=n}^{2n-1}\left| a_k - a_{k+3} \right|
	\geq \sum_{k \in C_n} \left| a_k - a_{k+3} \right|
	= \sum_{k \in C_n} \left| \frac{1}{k\ln(k+1)} - \frac{1}{k\ln(k+1)} - \frac{1}{(k+3)^{1+\frac{1}{p}}\ln(k+4)} \right|
	\end{equation*}
	\begin{equation*}
	= \sum_{k \in C_n} \frac{1}{(k+3)^{1+\frac{1}{p}}\ln(k+4)}
	\geq \frac{1}{(n+3)^{1+\frac{1}{p}}\ln(n+4)} \frac{n}{12}
	= \frac{1}{(n+3)^{1+\frac{1}{p}}\ln(n+4)} 
	\end{equation*}
	On the other hand, we get
	\begin{equation*}
	{}_3\beta(1) \leq C \frac{1}{n\ln(n+1)}
	\end{equation*}
	Therefore, the inequality
	\begin{equation*}
	\frac{1}{(n+3)^{1+\frac{1}{p}}\ln(n+4)} \leq C \frac{1}{n\ln(n+1)}
	\end{equation*}
	can not be satisfied if $n \rightarrow \infty$.
	
	Now, we will show that the series (\ref{sinus}) is divergent in $x_0 = \frac{2}{3}\pi$
	\begin{equation*}
	\sum_{k=1}^{6N+5} a_k \sin (kx_0) 
	= \sum_{k=1}^{N} \sum_{l=0}^{5} a_{6k+l} \sin((6k+l)\frac{2}{3}\pi)
	= \sum_{k=1}^{N} \left(  
	a_{6k} \sin(4\pi) 
	\right.
	\end{equation*}
	\begin{equation*}
	+ a_{6k+1} \sin((6k+1)\frac{2}{3}\pi)
	+ a_{6k+2} \sin((6k+2)\frac{2}{3}\pi) 
	+ a_{6k+3} \sin((6k+3)\frac{2}{3}\pi) 
	\end{equation*}
	\begin{equation*}
	\left.
	+ a_{6k+4} \sin((6k+4)\frac{2}{3}\pi)
	+ a_{6k+5} \sin((6k+5)\frac{2}{3}\pi)
	\right)
	\end{equation*}
	\begin{equation*}
	= \sum_{k=1}^{N} \left(
	a_{6k+1} \sin(4k\pi+\frac{2}{3}\pi)
	+ a_{6k+2} \sin(4k\pi+\frac{4}{3}\pi)
	\right.
	\end{equation*}
	\begin{equation*}
	\left.
	+ a_{6k+3} \sin(4k\pi+2\pi)
	+ a_{6k+4} \sin(4k\pi+\frac{8}{3}\pi)
	+ a_{6k+5} \sin(4k\pi+\frac{10}{3}\pi)
	\right)
	\end{equation*}
	\begin{equation*}
	= \sum_{k=1}^{N} \left(
	a_{6k+1} \sin(\frac{2}{3}\pi)
	+ a_{6k+2} \sin(\frac{4}{3}\pi)
	+ a_{6k+4} \sin(\frac{2}{3}\pi)
	+ a_{6k+5} \sin(\frac{4}{3}\pi)
	\right)
	\end{equation*}
	\begin{equation*}
	= \sum_{k=1}^{N} \left(
	a_{6k+1} \sin(\frac{2}{3}\pi)
	+ a_{6k+2} \left(-\sin(\frac{2}{3}\pi)\right)
	+ a_{6k+4} \sin(\frac{2}{3}\pi)
	+ a_{6k+5} \left(-\sin(\frac{2}{3}\pi)\right)
	\right)
	\end{equation*}
	\begin{equation*}
	= \sin(\frac{2}{3}\pi)\sum_{k=1}^{N} \left[ 
	\left( a_{6k+1} - a_{6k+2} \right) + \left( a_{6k+4} - a_{6k+5} \right)
	\right]
	\end{equation*}
	\begin{equation*}
	= \sin(\frac{2}{3}\pi)\sum_{k=1}^{N} \left[ 
	\left( \frac{3}{(6k+1)\ln(6k+2)} - \frac{1}{(6k+2)\ln(6k+3)} \right)
	\right.
	\end{equation*}
	\begin{equation*}
	\left.
	+ \left( \frac{3}{(6k+4)\ln(6k+5)} - \frac{1}{(6k+5)\ln(6k+6)} \right)
	\right]
	\end{equation*}
	\begin{equation*}
	\geq \sin(\frac{2}{3}\pi)\sum_{k=1}^{N} \left[ 
	\left( \frac{3}{(6k+2)\ln(6k+2)} - \frac{1}{(6k+2)\ln(6k+2)} \right)
	\right.
	\end{equation*}
	\begin{equation*}
	\left.
	+ \left( \frac{3}{(6k+5)\ln(6k+5)} - \frac{1}{(6k+5)\ln(6k+5)} \right)
	\right]
	\end{equation*}
	\begin{equation*}
	= \sin(\frac{2}{3}\pi)\sum_{k=1}^{N}
	\left( \frac{2}{(6k+2)\ln(6k+2)} + \frac{2}{(6k+5)\ln(6k+5)} \right)
	\end{equation*}
	\begin{equation*}
	\geq 4\sin(\frac{2}{3} \pi) \sum_{k=1}^{N} \frac{1}{(6k+5)\ln(6k+5)} \rightarrow \infty \text{ as }N\rightarrow \infty.
	\end{equation*} 
	This ends our proof. $\square$
	
\subsection{Proof of Theorem 11}
	The proof is similar to the proof of Theorem 10. The only difference is that we use \eqref{cosinus} and \eqref{cosinus_roznicowanie} instead of \eqref{sinus} and \eqref{sinus_roznicowanie}, respectively.

\subsection{Proof of Theorem 12}
	The proof is similar to the proof of Theorem 10. The only difference is that we use \eqref{exp} and Lemma \ref{roznicowanie_exp} instead of \eqref{sinus_roznicowanie} and Lemma \ref{sinus_roznicowanie}, respectively.

\end{document}